\documentclass[11pt]{article}%
\usepackage{amsfonts}
\usepackage{amsmath}
\usepackage{natbib}
\usepackage{graphicx}
\usepackage{here}
\usepackage[doublespacing]{setspace}
\usepackage{fullpage}%
\usepackage{amssymb}

\newcommand{\1}{1\hspace{-.55ex}\mbox{l}}
\numberwithin{equation}{section}
\begin{document}

\title{Using decomposed household food acquisitions as inputs of a Kinetic Dietary
Exposure Model.}
\author{Olivier Allais\thanks{INRA-CORELA, Laboratoire de recherche sur la
consommation, Ivry sur Seine, France.} and Jessica
Tressou\thanks{INRA-M\'{e}t@risk, M\'{e}thodologies d'analyse des risques
alimentaires, Paris, France and Hong Kong University of Science and
Technology, Department of Information and Systems Management, Hong Kong.}
\thanks{The second author research is in part supported by Hong Kong RGC Grant
\#601906.} \thanks{Corresponding author: Dr Jessica Tressou; Mail address:
Hong Kong University of Science and Technology, ISMT, Clear Water Bay,
Kowloon, Hong Kong; Email: \textit{tressou@inapg.fr.}}}
\date{}
\maketitle

\begin{abstract}
Foods naturally contain a number of contaminants that may have different and
long term toxic effects. This paper introduces a novel approach for the
assessment of such chronic food risk that integrates the pharmakokinetic
properties of a given contaminant. The estimation of such a Kinetic Dietary
Exposure Model (KDEM) should be based on long term consumption data which, for
the moment, can only be provided by Household Budget Surveys such as the
SECODIP panel in France. A semi parametric model is proposed to decompose a
series of household quantities into individual quantities which are then used
as inputs of the KDEM. As an illustration, the risk assessment related to the
presence of methylmercury in seafoods is revisited using this novel approach.


\end{abstract}

\newpage

\section*{Introduction}

The quantitative assessment of dietary exposure to certain contaminants is of
high priority to the Food and Agricultural Organization and the World Health
Organization (FAO/WHO). For example, excessive exposure to methylmercury, a
contaminant mainly found in fish and other seafood (mollusks and shellfish)
may have neurotoxic effects such as neuronal loss, ataxia, visual disturbance,
impaired hearing, and paralysis
\citep{WHO90}%
. Quantitative risk assessments for such chronic risk require the comparison
between a tolerable dose of the contaminant called Provisional Tolerable
Weekly Intake (PTWI) and the population's usual intake. The usual intake
distribution is generally estimated from independent individual food
consumption surveys (generally not exceeding 7 days) and food contamination
data. Several models have been developed to estimate the distribution of usual
dietary intake from short-term measurements
\citep[see for example, ][]{Nus,Hoff}%
.\emph{\ }The proportion of consumers whose usual weekly intake exceeds the
PTWI can then be viewed as a risk indicator
\citep[see for example, ][]{Mercure}%
. This kind of risk assessment does not account for the underlying dynamic
process, \textit{i.e.}\emph{\ }for the fact that the contaminant is ingested
over time and naturally eliminated at a certain rate by the human body.
Moreover, longer term measurements of consumption are available through
household budget surveys (HBS).

In this paper, we propose to use HBS data to quantify individual long term
exposure to a contaminant. This data provides long time series of household
food acquisitions which are first used in a decomposition model, similar to
the one proposed by \cite{Che97,Che98} in the nutrition field, in order to
obtain time series of individual intakes. Then, the pharmacokinetic properties
of the contaminant are integrated into an autoregressive model in which the
current body burden is defined as a fraction of the previous one plus the
current intake.

From a toxicological point of view, this approach is, to our knowledge, novel
and hence requires the definition of an ad-hoc long term safe dose as proposed
in the next section. We refer to this autoregressive model as Kinetic Dietary
Exposure Model (KDEM).

From a statistical point of view, such autoregressive models are well known in
general time series analysis
\citep[see for example, ][]{Ham94}
and most of the paper is devoted to the description of the decomposition
model. This statistical model aims at estimating individual quantities from
total household quantities and structures. This problem is similar to that
studied by \cite{Eng86}, \cite{Che97,Che98}, and \cite{VasTri00}, and is
addressed in a slightly different way. In the present article, the individual
contaminant intake is firstly viewed as a nonlinear function of age within
each gender, with time and socioeconomic characteristics being secondly
introduced in a linear way. The nonlinear function is represented by a
truncated polynomial spline of order 1 that admits a mixed model spline
representation
\citep[section 4.9 in ][]{RWC}%
. These choices yield a simple linear mixed model which is estimated by
REstricted Maximum Likelihood
\citep[REML, ][]{Pat71}%
. One major extension of the proposed model compared to \cite{Che97} is the
introduction of dependence between the individual intakes of a given household.

In the next section, focusing on the methylmercury example even though the
method is much more general and could be applied to any chronic food risk,
SECODIP data are described along with the construction of a household intake
series and the individual cumulative and long term exposure concepts yielding
the KDEM. Section \ref{methodo} is devoted to the statistical methodology used
to decompose the household intake series into individual intake series, namely
the presentation of the model and its estimation and tests. Section
\ref{results} displays the results for the quantification of long term
exposure to methylmercury of the French population using the 2001 SECODIP
panel. Finally, a discussion on the use of household acquisition data, with
the focus on the French SECODIP panel, is conducted in section \ref{discus}
with respect to the proposed long term risk analysis.

\section{Motivating example: risk related to methylmercury in seafoods in the
French population\label{LT_data}}

In this section, the Kinetic Dietary Exposure Model (KDEM) and the concept of
long term risk are defined. Then a brief panorama of consumption data in
France is given and the way the SECODIP HBS data will be used as an input of
the KDEM is described.

\subsection{Cumulative exposure and long term risk: the Kinetic Dietary
Exposure Model (KDEM)}

The main objective of the analysis is to assess individuals' long term
exposure to a contaminant to deduce whether these individuals are at risk or
not. As mentioned in the introduction the only \textquotedblright safe
dose\textquotedblright\ reference is the PTWI expressed in terms of body
weight (\textit{relative} intake). Unfortunately, TNS SECODIP did not record
the body weight of the individuals until 2001. The body weights are thus
estimated from independent data sets; namely the French national survey on
individual consumption
\citep[INCA, ][]{INCA}
for people older than 18, and the weekly body weight distribution available
from French health records (\cite{SemPed}) for individuals under 18. In both
cases, gender differentiation is introduced.

Assume that estimations of the individual weekly intakes are available, that
is $y_{i,h,t}$ denotes the intake of individual $i$ belonging to household $h
$ for the $t^{th}$\ week (with $i=1,\emph{\ldots},n_{h,t};$ $h=1,\emph{\ldots
},H$ and $t=1,\emph{\ldots},T),$ and $D_{i,h,t}$\ denotes the same quantity
expressed on a body weight basis. The cumulative exposure up to the $t^{th}%
$\ week of this individual is then given by
\begin{equation}
S_{i,h,t}=\exp(-\eta)\cdot S_{i,h,t-1}+D_{i,h,t},\label{Cumul_expo}%
\end{equation}
where $\eta>0$ is the natural dissipation rate of the contaminant in the
organism. This dissipation parameter is defined from the so called
\textit{half life} of the contaminant,which is the time required for the body
burden to decrease by half in the absence of any new intake. For
methylmercury, the half life, denoted by $l_{1/2},$ is estimated to $6$ weeks,
so that $\eta=\ln(2)/l_{1/2}:=\ln(2)/6$
\citep{elimMeHg}%
.

The autoregressive model defined by $\left(  \ref{Cumul_expo}\right)  $ and a
given initial state $S_{i,h,0}=D_{i,h,0}$ has a stationary solution since
$\exp(-\eta)<1.$ As a convention, $S_{i,h,0}$ is set to the mean of all
positive exposures $\left(  D_{i,h,t}\right)  _{t=1,\emph{\ldots},T}$.
However, this convention has little impact on the level of an individual's
long term exposure since the contribution of the initial state $S_{i,h,0}$
tends to zero as $t$ increases. We call this autoregressive model "KDEM" for
Kinetic Dietary Exposure Model.

The individual cumulative exposure $S_{i,h,t}$ can be considered to be the
long term exposure of an individual for sufficiently large values of $t$. For
methylmercury, the long term steady state of the individual exposure to a
contaminant is reached after $5$ or $6$ half lives according to Dr P.
Granjean, a methylmercury expert. Thus, the long term individual's exposure to
methylmercury is defined as the cumulative exposure reached after say
$6l_{1/2}=36$ weeks.

The risk assessment usually consists of comparing the exposure with the so
called Provisional Tolerable Weekly Intake (PTWI). This tolerable dose,
determined from animal experiments and extrapolated to humans, refers to the
dose an individual can ingest throughout his entire life without appreciable
risk. For methylmercury, the PTWI is set to $1.6$ microgram per kilogram of
body weight per week
\citep[$1.6$ $\mu$g/kg bw, see ][]{JECFAMeHg}%
.

In our dynamic approach, the long term exposure is compared to a reference
long term exposure denoted by $S^{ref}$, and defined as the cumulative
exposure of an individual whose weekly intake is equal to the PTWI, $d$, such
as%
\begin{equation}
S^{ref}=\lim_{t\rightarrow\infty}S_{t}^{ref}=\frac{d}{1-\exp(-\eta
)},\label{Cumulref_lim}%
\end{equation}
where%
\begin{equation}
S_{t}^{ref}=\sum_{s=0}^{t}d\exp(-\eta(t-s))=d\frac{\exp(-\eta(t+1))-1}%
{\exp(-\eta)-1}.\label{cumul_ref}%
\end{equation}
For methylmercury, the reference for long term exposure $S^{ref}$ is $14.6$
$\mu$g/kg bw. An individual is then assumed to be at risk if his cumulative
exposure $S_{i,h,t}$ exceeds the reference $S_{t}^{ref}$ for any $t>6l_{1/2}$.

This KDEM model requires some long surveys of individual intakes which are not
monitored and can only be approximated from available consumption data and
contamination data.

\subsection{From household acquisition data to household intake series}

Two current major consumption data sources in France are the national survey
on individual consumption
\citep[INCA, ][]{INCA}
and the SECODIP panel managed by the company TNS SECODIP. Most quantitative
risk assessments conducted by the French agency for food safety (AFSSA) use
the 7 day individual consumption data of the INCA survey jointly with
contamination data collected by several French institutions. Regarding
methylmercury, seafood contamination data have been collected through
different analytical surveys
(\citeauthor{MAAPAR}, \citeyear{MAAPAR}; \citeauthor{IFREMER}, \citeyear
{IFREMER})
and were used in \cite{Mercure} and \cite{MercureAC} combined with the INCA
survey. In this paper, a methodology using the SECODIP data is developed
\citep[see ][ for a full description of this database]{Boi05}%
.

The company TNS SECODIP has been collecting the weekly food acquisition data
of about five thousand households since 1989. All participating households
register grocery purchases through the use of EAN bar codes but other grocery
purchases are registered differently: the fresh fruit and vegetable purchases
are recorded by the FL sub-panel while fresh meat, fresh fish and wine
purchases are recorded by the VP sub-panel. The households are selected by
stratification according to several socioeconomic variables and stay in the
survey for about 4 years. TNS SECODIP provides weights for each sub-panel and
each period of 4 weeks to make sure of the representativeness of the results
in terms of several socioeconomic variables. TNS SECODIP also defines the
notion of household activity which refers to the correct and regular reporting
of household purchases over a year. For each household, the age and gender of
each member of the household are retained in our decomposition model with some
socioeconomic variables: the region, the social class (from modest to
well-to-do), the occupation category and level of education of the principal
household earner.

For methylmercury risk assessment, the households of the VP panel are
considered; in the 2001 data set, there are $H=3229$ active households
(corresponding to $9288$ individuals) and $T=53$ weeks during which the
households may or may not acquire seafood. The weekly purchases of seafood are
clustered into two categories (\textquotedblright Fish\textquotedblright\ and
\textquotedblright Mollusks and Shellfish\textquotedblright) for which the
mean contamination levels are calculated from the MAAPAR-IFREMER data and are
given in table \ref{conta}.

\begin{center}
\emph{Table \ref{conta} around here, see page \pageref{conta}}
\end{center}

Household intake series ($\left(  y_{h,t}\right)  _{h=1,...,H;t=1,...,T}$) are
computed as the cross product between weekly purchases of seafoods which are
assimilated to weekly consumptions, and mean contamination levels. They are
expressed in micrograms per week ($\mu g/w$). The food "purchase-consumption"
assimilation is of course arguable and will be the main subject of the final
discussion (see section \ref{discus}). An additional assumption concerns the
household size, denoted by $n_{h,t}$ for the household $h$ and the week $t.$
This can indeed vary over time in the case of a birth or death of a household
member. Since a new born baby will not consume fish in his first few months,
we assume that food diversification (and hence consumption of seafoods) starts
at one year of age, yielding a total sample of $8913$ individuals for the 2001
panel. These household intake series are then decomposed into individual
intake series using the model described in the next section. These individual
intake series are then used as imputs of the KDEM.

\section{Statistical methodology\label{methodo}}

In this section, the decomposition model is described and compared to similar
models described in the literature, namely \cite{Che97,Che98,VasTri00}. Its
estimation and some structure tests are then presented.

\subsection{The decomposition model}

\subsubsection{General principle}

Consider a household composed of $n_{h,t}$ members, each member having
unobserved weekly intakes $y_{i,h,t}$, with $i=1,$\emph{\ldots}$,n_{h,t},$
$h=1,$\emph{\ldots}$,H$, and $t=1,$\emph{\ldots}$,T$. The week $t$ intake of a
household $h$ is simply the sum across household members of the individual
weekly intakes, such as
\begin{equation}
y_{h,t}=\sum_{i=1}^{n_{h,t}}y_{i,h,t}.\label{sum}%
\end{equation}
As detailed below, the individual weekly intake $y_{i,h,t}$ is assumed to
depend on

\begin{itemize}
\item the age and gender of the individual via a function $f,$

\item some socioeconomic characteristics of the household,

\item time (seasonal variations).
\end{itemize}

There are obviously several ways to model the individual intake under these
assumptions and this choice leads to more or less simple estimation
procedures. In \cite{Che97,Che98,VasTri00}, a discretization argument on age
is used leading to a penalized least square estimation of a great number of
parameters, that is one parameter for each year of age and gender. We propose
to use a truncated polynomial spline of order 1 for each gender, which admits
a mixed model spline representation for $f.$ As far as socioeconomic
characteristics are concerned, \cite{Che97} retained a multiplicative
specification whereas \cite{VasTri00} chose the additive one. In the
multiplicative model, a change in income for example would proportionally
affect all the individual intakes whereas in the additive setting, they would
be affected by the same value. Following \cite{VasTri00}, we retained the
additive specification since the difference between the two specifications may
not be notable, and the additive setting yields to a much simpler estimation
procedure (linear model). Finally, time dependency is only introduced in
\cite{Che98} to track changes with age within cohorts: this time dependency is
directly introduced into the function $f$ that is bivariately smoothed
according to age and time
\citep[cf. ][]{Gre94}%
. Again, we adopt a simpler specification in which time is introduced as a
dummy variable. All these assumptions yield an individual model of the form%
\begin{equation}
y_{i,h,t}=x_{i,h,t}\beta+z_{i,h,t}u+w_{h,t}\gamma+\delta_{t}\alpha
+\varepsilon_{i,h,t},\label{modelI}%
\end{equation}
where the terms $x_{i,h,t}\beta+z_{i,h,t}u$ stand for the mixed model spline
representation of the function $f,$ the term $w_{h,t}\gamma$ denotes the
socioeconomic effects, the term $\delta_{t}\alpha$ the time effect, and
$\varepsilon_{i,h,t}$ is the individual error term.

Combining $\left(  \ref{sum}\right)  $ and $\left(  \ref{modelI}\right)  ,$ we
obtain the final rescaled household model given by%
\begin{equation}
Y_{h,t}=X_{h,t}\beta+Z_{h,t}u+\sqrt{n_{h,t}}w_{h,t}\gamma+\sqrt{n_{h,t}}%
\delta_{t}\alpha+\varepsilon_{h,t},\label{HousMod}%
\end{equation}
where $Y_{h,t}\equiv\sum_{i=1}^{n_{h,t}}y_{i,h,t}/\sqrt{n_{h,t}},$
$X_{h,t}\equiv\sum_{i=1}^{n_{h,t}}x_{i,h,t}/\sqrt{n_{h,t}},$ $Z_{h,t}%
\equiv\sum_{i=1}^{n_{h,t}}z_{i,h,t}/\sqrt{n_{h,t}},$ and $\varepsilon
_{h,t}\equiv\sum_{i=1}^{n_{h,t}}\varepsilon_{i,h,t}/\sqrt{n_{h,t}}.$

\subsubsection{Specification details}

\paragraph{Age-gender function specification}

Let $a_{i,h,t}$ and $s_{i,h}$ denote the age and sex of individual $i$ of
household $h$ for the $t^{th}$ week. Individual dietary intake is generally
different according to the gender of individuals, so the function $f$ takes
the following form%
\[
f(a_{i,h,t},s_{i,h})=f_{M}(a_{i,h,t})%
\1%
_{\left\{  s_{i,h}=M\right\}  }+f_{F}(a_{i,h,t})%
\1%
_{\left\{  s_{i,h}=F\right\}  },
\]
where $f_{M}(.)$ and $f_{F}(.)$ are age-intake relationships for males (M) and
females (F) respectively, and $%
\1%
_{\{A\}}$ is the indicator function of event $A$. The function $f_{S}(.)$ is
approximated by a spline of order one with a truncated polynomial basis for
either sex, such as%
\begin{equation}
f_{S}(a_{i,h,t})=\beta_{0}^{S}+\beta_{1}^{S}a_{i,h,t}+\sum_{k=1}^{K_{S}}%
u_{k}^{S}\left(  a_{i,h,t}-\kappa_{S,k}\right)  _{+},\label{deffS}%
\end{equation}
where the $\left(  \kappa_{S,k}\right)  _{k=1,\emph{\ldots},K_{S}}$ are nodes
chosen from an age list and
\[
\left(  a_{i,h,t}-\kappa_{S,k}\right)  _{+}\equiv\left(  a_{i,h,t}%
-\kappa_{S,k}\right)
\1%
_{\left\{  a_{i,h,t}-\kappa_{S,k}>0\right\}  }%
\]
denotes the positive part of the difference between the age of the individual
$a_{i,h,t}$ and the node $\kappa_{S,k}$ and the $u_{k}^{S}$ are random effects
assumed to be i.i.d. Gaussian with distribution$\ \mathcal{N}\left(
0,\sigma_{u_{S}}^{2}\right)  $. This last assumption allows us to introduce
some penalties into the model and to smooth the function $f_{S}$ yielding a
mixed model representation for the spline as shown in
\cite{Spe91,Ver99,Bru99,RWC}. As in \cite{RWC}, page 125, the total number of
nodes $K_{S}$ is set to $\min\left\{  \left\vert \frac{a_{S,d}}{4}\right\vert
,35\right\}  $, where $a_{S,d}$ is the list of distinct ages for individuals
of sex $S$, and the nodes $\kappa_{S,k}$ are defined as the $\left(
\frac{k+1}{K_{S}+2}\right)  ^{th}$ percentile of vector $a_{S,d}$ for
$k=1,$\emph{\ldots}$,K_{S}$.

Defining $x_{i,h,t}$ as a line vector $\left(
\begin{array}
[c]{cccc}%
\1%
_{\left\{  s_{i,h}=M\right\}  } & a_{i,h,t}%
\1%
_{\left\{  s_{i,h}=M\right\}  } &
\1%
_{\left\{  s_{i,h}=F\right\}  } & a_{i,h,t}%
\1%
_{\left\{  s_{i,h}=F\right\}  }%
\end{array}
\right)  ,$ and $z_{i,h,t}$ as the line vector $\left\{  \left(
a_{i,h,t}-\kappa_{S,k}\right)  _{+}%
\1%
_{\left\{  s_{i,h}=S\right\}  }\right\}  _{k=1,\emph{\ldots},K_{S}\text{
};\text{ }S=M,F},$ we finally obtain the first terms of $\left(
\ref{modelI}\right)  ,$ that is $f(a_{i,h,t},s_{i,h})=x_{i,h,t}\beta
+z_{i,h,t}u.$

\paragraph{Socioeconomic characteristics and time dependency}

In the application, all the socioeconomic characterics are categorical
variables. Consider the $Q$ categorical variables $W_{h,t}^{(q)},$
$q=1,\ldots,Q,$ with $m_{q}$ modalities, and fix the $m_{q}^{th}$ modality as
the reference modality, then the socioeconomic effect term in $\left(
\ref{modelI}\right)  $ and $\left(  \ref{HousMod}\right)  $ is%
\[
w_{h,t}\gamma=\sum_{q=1}^{Q}\sum_{m=1}^{m_{q}-1}\gamma_{q,m}%
\1%
_{\left\{  W_{h,t}^{(q)}=m\right\}  ,}%
\]
where $\gamma_{q,m}$ is the effect of the $m^{th}$ modality of the
socioeconomic variable $q.$

Similarly, time is only measured by weekly counts throughout the year so that
the time effect in $\left(  \ref{modelI}\right)  $ and $\left(  \ref{HousMod}%
\right)  $ is simply%
\[
\delta_{t}\alpha=\sum_{\substack{\tau=1 \\\tau\neq\tau_{R}}}^{T}\alpha_{\tau}%
\1%
_{\left\{  \tau=t\right\}  },
\]
where $\alpha_{\tau}$ is the effect of week $\tau$ and $\tau_{R}$ is the
reference week.

\paragraph{Error specification}

The error at the individual level $\varepsilon_{i,h,t}$ is assumed to be
Gaussian with zero mean, and the variance-covariance structure is such that

\begin{itemize}
\item households are independent, i.e. $\forall i,i^{\prime},t,t^{\prime}%
$\ and $\forall h\neq h^{\prime}$%
\[
cov(\varepsilon_{i,h,t},\varepsilon_{i^{\prime},h^{\prime},t^{\prime}})=0,
\]

\item members of the same household are dependent, that is for $\forall h,t$
and $i\neq i^{\prime},$%
\[
cov(\varepsilon_{i,h,t},\varepsilon_{i^{\prime},h,t})=\rho\sigma_{\varepsilon
}^{2},
\]
where $\rho$ measures the dependence between individuals within the same household.

\item there is no time dependence, that is $\forall i,i^{\prime}$\ and
$\forall t\neq t^{\prime}$
\[
cov(\varepsilon_{i,h,t},\varepsilon_{i^{\prime},h,t^{\prime}})=0.
\]

\end{itemize}

In the rescaled household model $\left(  \ref{HousMod}\right)  $, the error
$\varepsilon_{h,t}\equiv\sum_{i=1}^{n_{h,t}}\varepsilon_{i,h,t}/\sqrt{n_{h,t}%
}$ is i.i.d. Gaussian with a zero mean and a variance $R$ such that $\forall
t,t^{\prime}$ and $\forall h\neq h^{\prime}$,%
\begin{equation}
\mathbb{V}(\varepsilon_{h,t})=\rho\sigma_{\varepsilon}^{2}n_{h,t}%
+(1-\rho)\sigma_{\varepsilon}^{2}\text{ and }cov(\varepsilon_{h,t}%
,\varepsilon_{h^{\prime},t^{\prime}})=0.\label{Var_cov hous}%
\end{equation}

\subsection{Estimation and tests\label{estim}}

The model $(\ref{HousMod})$ is a linear mixed model that can be estimated
using restricted maximum likelihood (REML) techniques, see \cite{RWC} for
details. An attractive consequence of the use of the mixed model
representation of a penalized spline in $\left(  \ref{deffS}\right)  $ is that
mixed model methodology and software can be used to estimate the parameters
and predict the random effect in the resulting household model. The amount of
smoothing of the underlying functions $f_{S}$ is estimated with the REML
technique via the estimation of $\sigma_{u_{S}}^{2}$. The estimation was
conducted using \textregistered SAS \texttt{MIXED} procedure. To get
estimators for $\sigma_{\varepsilon}^{2}$ and $\rho,$ asymptotic least square
techniques combined with the linear relationship between the variance given in
$\left(  \ref{Var_cov hous}\right)  $ and the household size were used. More
precisely, a residual variance $\sigma_{n}^{2}$ is first estimated for each
household size $n=1,\ldots,N=\max n_{h,t}$ using an option of the
\texttt{MIXED} procedure (see the program for the detailed syntax). Then,
ordinary least square regression and the delta method give estimators for
$\sigma_{\varepsilon}^{2}$ and $\rho$ and their standard deviations.

The individual intake is then predicted by
\begin{equation}
\widehat{y_{i,h,t}}=x_{i,h,t}\widehat{\beta}+z_{i,h,t}\widehat{u}+w_{h,t}%
\hat{\gamma}+\delta_{t}\hat{\alpha},\label{EstIndExpo}%
\end{equation}
where $\widehat{\beta}$, $\hat{\gamma}$, and $\hat{\alpha}$ are the estimators
of $\beta$, $\gamma$, and $\alpha$ respectively and $\widehat{u}$ is the best
prediction of the random effect $u$ in the model $\left(  \ref{HousMod}%
\right)  $.

Confidence and prediction intervals can be built for the prediction
$\widehat{y_{i,h,t}}$ as proposed in \cite{RWC} and several tests can be
conducted in this model:

\begin{enumerate}
\item Are the random effects different according to sex? In other words, is
the assertion $\sigma_{u_{M}}^{2}=\sigma_{u_{F}}^{2}=\sigma_{u}^{2}$ true?

\item Another question is the necessity for such random effects. Is the
assertion $\sigma_{u}^{2}=0$ (resp. $\sigma_{u_{M}}^{2}=0$ or $\sigma_{u_{F}%
}^{2}=0)$ true?

\item More globally, is the function $f$ the same for both sexes? Is the
assertion $f_{M}=f_{S}$ true?
\end{enumerate}

These tests can be conducted using classical likelihood (or restricted
likelihood) ratio techniques. The likelihood ratio statistic is asymptotically
distributed as a chi square with a degree of freedom being the number of
tested equalities, except for point 2, where the limiting distribution is
known to be a mixture of chi-square
\citep{seli87,crarup03}
because the test concerns the frontier of the parameter definition
($\sigma_{u}^{2}\in\left[  0,+\infty\right[  $).

\section{Applying our methodology to the methylmercury risk
assessment\label{results}}

In this section, we illustrate our approach on our motivating example.
Firstly, several tests are conducted on the decomposition model, and secondly,
individual long term exposure is compared to the reference long term exposure
described in section \ref{LT_data}.

\subsection{Estimation and tests on the structure of the model}

Table \ref{Esti_gen_model} shows the REML estimates for all socioeconomic
variables (parameter $\gamma)$ and the p-values of Student tests in the model
$(\ref{HousMod})$. The socioeconomic variables used are household income,
region of residence, occupation category and level of education of the
principal household earner. For each socioeconomic variable, the reference
modality is given in Table \ref{Esti_gen_model}. We assume here that

\begin{itemize}
\item the function $f$ differs according to the gender but the random effect
does not ($f_{M}\neq f_{F}$ and $\sigma_{u_{M}}^{2}=\sigma_{u_{F}}^{2}),$

\item the maximum household size $\overline{N}$ is set to $6$ for
variance-covariance estimation. Indeed, the dependence between individuals
within the same household depends on the household size $n_{h}$\ in
$(\ref{Var_cov hous})$. For each household size, a variance is estimated, and
estimates of $\rho$\ and $\sigma^{2}$\ are obtained using asymptotic least
square techniques as mentioned in section \ref{estim}. Since large households
are not numerous in the database, the estimations are implemented with a
maximum household size, $\overline{N}$, set to $6$; it is assumed that there
is a common variance for all households with size greater than $\overline{N}$.
\end{itemize}

In this sub-section, we show the results of several tests we carried out to
simplify the interpretation of our study. These tests have been implemented in
a hierarchical way, starting with the highest-order interaction terms,
combining to the reference modality the modality which does not differ
significantly from the reference. All tests are performed on the $5\%$ level
of significance and each new hypothesis is tested, conditionally on the
results of the previous tests. Each null hypothesis and the p-value resulting
from the appropriate F-test are shown in Table \ref{Test_socio}.

First of all, concerning the occupation category variable, the self-employed
modality does not significantly differ from the reference modality blue collar
workers $(H1,$ $Pval=0.771)$. Refitting the model with the reference modality
\textquotedblright Blue collar workers and self employed\textquotedblright,
all the socioeconomic variables are significantly different from the
reference. Then, F-tests allow us to conclude that the resulting three groups
are significantly different from each other $(H2,$ $H3,$ $H4)$.

Let us now consider the region of residence variable. First, there are some
very substantial differences among the $4$ regions of residence $(H5,$
$Pval=<0.001)$. However, the modality "North, Brittany, and Vendee coast" and
the modality "Paris and its suburbs" should be grouped $(H6$ $c$,
$Pval_{c}=0.881)$. Then, the other tests implemented for the level of
education and income variables suggest that no further simplification is
possible (see p-values of null hypotheses $H7$, $H8$, $H9$ in Table
\ref{Test_socio}). Finally, the overall F-test comparing our resulting final
model to the original model $(\ref{HousMod})$ shows that no important variable
has been left out of the model $\left(  Pval=0.59\right)  $.

Table \ref{Esti_final_model} shows the parameter estimates and p-values of the
Student's t-tests for all socioeconomic variables of the reduced final model.
The income effects on individual exposure are those expected: the richer the
households are, the higher their exposures are because seafoods are expensive.
Furthermore, living in a coastal region or in Paris and its suburbs brings
about larger individual exposure relatively to living in a non coastal region
because of the more ready supply of seafoods in these regions. Moreover, the
more educated you are, the larger the individual exposure is. The occupation
category of the principal household earner has an unexpected effect on the
individual exposure. Indeed a higher exposure is expected for white collar
workers and retirees whan compared to blue collar workers but an opposite
effect is observed. This may be explained by the fact that the reference
modality for this variable is a very heterogeneous modality also comprising
managers and self-employed persons (farmers and craftsmen). Another
explanation could be that white collars workers have a higher propensity to
eat out in restaurants whereas outside the home consumption is not included in
the model.

\begin{center}
\emph{Table \ref{Esti_gen_model} around here, see page
\pageref{Esti_gen_model}}

\emph{Table \ref{Test_socio} around here, see page \pageref{Test_socio}}

\emph{Table \ref{Esti_final_model} around here, see page
\pageref{Esti_final_model}}
\end{center}

Likelihood ratio tests are implemented to test the structure of the final
model. First, the dependence of individual exposures to methylmercury within a
household is tested. The null hypothesis $\rho=0$ (cf. equation $\left(
\ref{Var_cov hous}\right)  $) is rejected (null Pval) which confirms that
individuals within the same household have correlated exposures. Then, we test
if the function $f$ is the same for both genders. The null hypothesis
$f_{M}=f_{F}$ is rejected (null Pval) but the null hypothesis $\sigma_{u_{M}%
}^{2}=\sigma_{u_{F}}^{2}$ is accepted. This means that individual exposure
differs with gender but both functions need the same amount of smoothing.

\subsection{The cumulative and the long term individual exposure}

\begin{figure}
[h]
\begin{center}
\includegraphics[
natheight=6.246500in,
natwidth=10.085400in,
height=3.8424in,
width=6.1869in
]%
{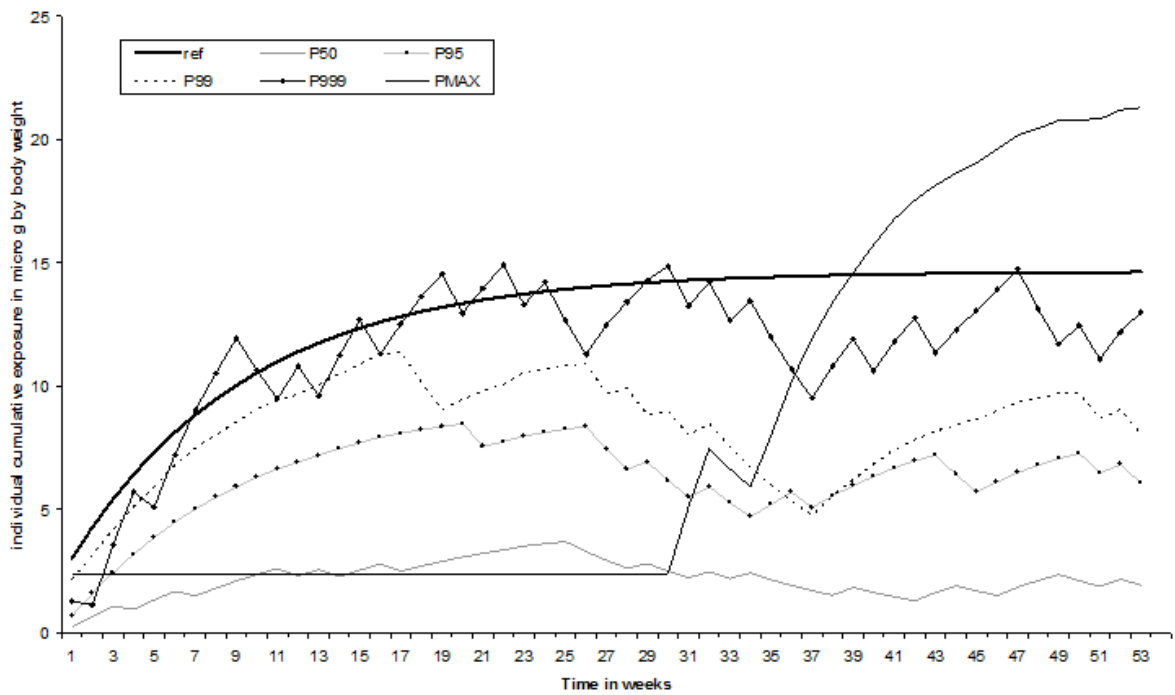}%
\caption{Cumulative exposure to MeHg (unit: $\mu$g per kg of body weight)}%
\label{Dist_expocum}%
\end{center}
\end{figure}

The cumulative individual exposure $S_{i,h,t}$ is calculated from the
estimated individual weekly intakes according to equation $(\ref{Cumul_expo})
$ and the resulting values for $t>35$ are compared to the reference cumulative
exposure defined by $(\ref{cumul_ref})$. Figure \ref{Dist_expocum} shows the
cumulative individual exposure over the $53$ weeks of the year $2001$ for
different individuals. Only certain percentiles of the distribution of the
individual cumulative exposures of the last week are displayed. For example,
the curve \texttt{Pmax} represents the cumulative exposure of an individual
whose last week's cumulative exposure is the highest. This is the cumulative
exposure of a girl who turned one year old during the 30th week of 2001, lives
in Paris or its suburbs in a well to do household.\emph{\ }

Very few individuals have a cumulative individual exposure above the reference
long term exposure. We estimate that only $0.186\%$ of individuals are deemed
at risk.\emph{\ }This risk index should be compared to the more common one
defined as the percentage of weekly intakes $D_{i,h,t}$\ exceeding the PTWI,
denoted $R_{1.6}$, such as $R_{1.6}=\frac{1}{nT}\sum_{t=1}^{T}\sum_{h=1}%
^{H}\sum_{i=1}^{n_{h}}%
\1%
\left(  D_{i,h,t}>1.6\right)  $. $R_{1.6}$\ is equal to $0.45\%$, and is
slightly higher since each occasional deviation above the PTWI increases the
risk index whereas only long term deviations above this PTWI should be taken
into account to assess the risk.

A deeper analysis of at risk individuals shows that all these vulnerable
individuals are children less than three years old. They represent $5.29\%$ of
the children aged between 1 and 3 in 2001. Further, no child of a modest
households is found to be at risk.

\begin{center}
\emph{Figure \ref{Dist_expocum} around here, see page \pageref{Dist_expocum}}
\end{center}

\section{Discussion\label{discus}}

As mentioned in section \ref{LT_data}, the use of household acquisition data
in a food safety context, and in our case the use of the SECODIP database for
assessing methylmercury dietary intakes, gives rise to some approximations:

\begin{enumerate}
\item Consumption outside of the home is out of the scope of household
acquisition data. TNS SECODIP does not provide any information on the
quantities of seafoods consumed out of the home or bought for outside
consumption. Nevertheless, \cite{Ser03} assert that these data are good
estimates for the consumption of the whole household. \cite{VasTri00} avoid
this question by using the term "availaibility" instead of intake or
consumption. However, as in \cite{Che97}, auxiliary information about outdoor
consumption could be introduced in the model as a correction factor accounting
for the propensity to eat outside of the home according to age, sex or
socioeconomic variables. The French INCA survey on individual consumptions
gives details about inside / outside the home consumption for $3003$
individuals people aged 3 and older. The mean outside the home consumption
proportion is $20\%$ for seafoods. Applying such a factor to all household
intakes yields a long term risk of $0.226\%$, and $R_{1.6}=0.791\%$.
Furthermore, in this case, a small proportion of consumers older than 3 years
old are vulnerable. Nevertheless, children aged between 1 and 3 in 2001 still
represent the most vulnerable consumer group, at $10\%$\ of the corresponding population.

\item The amount of food bought by a household can be different from the
amount actually consumed. Indeed, namely for seafoods, a non negligible part
is not edible: \cite{Ciqual} show than on average only 61\% of fresh or frozen
fish is edible. Besides, \cite{Dec94} also demonstrate some part of the
purchased food is thrown away, which also reduces the actual amount of food
consumed by a household. However, SECODIP does not specify whether the
quantity of fresh or frozen fish bought is ready to be consumed or as a whole
fish that needs some preparation. Applying such a factor to all household
intakes yields a long term risk of $0$.$00\%$, and $R_{1.6}=0.043\% $. If both
the $20\%$ outside of the home consumption correction factor and the $61\%$
edible proportion factor are applied to our series, the long term risk is
equal to $0.021\%$, $R_{1.6}=0.13\%$, and $1.06\%$ of the population of
children aged between 1 and 3 are vulnerable. These results stress that
applying such a correction factor to assess the actual quantity consumed is
probably too strong and is certainly a crude approximation of the quantity of
seafoods ingested. Thus, a more detailed database on fish and seafood is
needed, to realize an accurate assessment of exposure to methylmercury, taking
into account only the edible part of fish and other seafood.

Body weight information is crucial in a food safety context and will be
included in the future SECODIP data since it has now been added to the list of
required individual characteristics. The measurement error afferent to this
quantity will remain however, namely for children whose body weight changes a
lot throughout a year. Nevertheless, approximating the weekly body weight of
young children by the median of the weekly body weight distribution available
in French health records is the best approximation possible.

\item The food nomenclature of the SECODIP database is not as detailed as the
contamination database. Unfortunately, fish and seafood species are not well
documented so it is not possible to consider more than two food categories
when computing household intakes. This problem of nomenclature matching is
ubiquitous of food risk assessments since contamination analysis are generally
conducted independently from the food nomenclature of consumption data.
\end{enumerate}

These arguments mainly show the disadvantages of the use of household food
acquisition data such as the SECODIP database. Nevertheless, they also present
many advantages compared to the individual food record survey mainly used in
France in the food safety context:

\begin{itemize}
\item As mentioned before, households respond for a long period of time (the
average is 4 years in the SECODIP panel) which allows us to observe long term
behaviors and avoid some well known biases of individual food record surveys.
For example, respondents might over- (under-) declare certain foods with a
good (bad) nutritional value either deliberately or just because they
increased (reduced) their consumption for the short (7 days) period of the survey.

\item The individual surveys are expensive and very difficult to conduct.
Highly trained interviewers are required and extraordinary cooperation is
required from respondents. Household food acquisition data can serve many
other applications (economics or marketing) and, at least for the SECODIP
data, acquisition recording is simplified by optical scanning of food barcodes.
\end{itemize}

\section*{Conclusion}

In this paper, we proposed a methodology to assess chronic risks related to
food contamination using the example of methylmercury exposure through seafood
consumption. This methodology includes the definition of a Kinetic Dietary
Exposure Model (KDEM) that integrates the fact that contaminants are
eliminated from the body at different rates, the rate being measured by the
half life of the contaminant. In this paper, the estimation is based on the
use of household food acquisition data which are first decomposed into
individual intake data through a disaggregation model accounting for the
dependence among household members. Several extensions of this methodology are
currently studied. First, the disaggregation model could be improved by
considering a preliminary step in which we determine what member is an actual
consumer, in the spirit of the Tobit model. The KDEM idea is also currently
being developed by studying the stability and ergodic properties of the
underlying continuous time piecewise deterministic Markov process
\citep{BerClemTres}%
. The parameters of this new model are the intake distribution, the inter
intake time distribution and the dissipation rate distribution. In this
framework, the dissipation parameter $\eta$ of the KDEM model is random and
the intake and inter-intake distributions can be estimated either from
individual (INCA-type) data or household (SECODIP-type) data.

\bibliographystyle{ifac}
\bibliography{pra}
\newpage

\begin{center}
{\LARGE Figures and Tables}

\bigskip%

\begin{table}[h] \centering
\caption{Description of the contamination database (Unit:
microgram per kilogram}%
\begin{tabular}
[c]{|r|c|c|c|c|c|}\hline
& {\scriptsize Mean} & {\scriptsize Min} & {\scriptsize Max} &
{\scriptsize Standard Deviation} & {\scriptsize Number of analysis}\\\hline
{\scriptsize Fish} & {\scriptsize 0.147} & {\scriptsize 0.003} &
{\scriptsize 3.520} & {\scriptsize 0.235} & {\scriptsize 1350}\\\hline
{\scriptsize Mollusk and Shellfish} & {\scriptsize 0.014} &
{\scriptsize 0.001} & {\scriptsize 0.172} & {\scriptsize 0.011} &
{\scriptsize 1293}\\\hline
\end{tabular}
\label{conta}%
\end{table}%

\end{center}%

\begin{table}[h] \centering
\caption
{Restricted maximum likelihood estimates (REML) for age and all socioeconomic variables and the p-value of the Student's tests (Pval)\label
{Esti_gen_model}}%
\begin{tabular}
[c]{llll}\hline\hline
{\scriptsize Effect} & {\scriptsize Parameter} & {\scriptsize REML} &
{\scriptsize Pval}\\
{\scriptsize Income} & \multicolumn{3}{l}{{\scriptsize \ (ref: Mean sup)}}\\
{\scriptsize \quad Well to do} & ${\scriptsize \gamma}_{1}$ &
{\scriptsize 6.027} & {\scriptsize
$<$%
0.001}\\
{\scriptsize \quad Mean inf} & ${\scriptsize \gamma}_{2}$ &
{\scriptsize 2.686} & {\scriptsize
$<$%
0.001}\\
{\scriptsize \quad Modest} & ${\scriptsize \gamma}_{3}$ & {\scriptsize -1.928}
& {\scriptsize
$<$%
0.001}\\
{\scriptsize Region of residence } & \multicolumn{3}{l}{{\scriptsize (ref:
Noncoastal regions)}}\\
{\scriptsize \quad North, Brittany, Vendee coast } & ${\scriptsize \gamma}%
_{4}$ & {\scriptsize 0.962} & {\scriptsize 0.003}\\
{\scriptsize \quad South West coast} & ${\scriptsize \gamma}_{5}$ &
{\scriptsize 5.232} & {\scriptsize
$<$%
0.001}\\
{\scriptsize \quad Mediterranean coast} & ${\scriptsize \gamma}_{6}$ &
{\scriptsize 2.303} & {\scriptsize
$<$%
0.001}\\
{\scriptsize \quad Paris and its suburbs} & ${\scriptsize \gamma}_{7}$ &
{\scriptsize 1.023} & {\scriptsize 0.009}\\
{\scriptsize Occupation category of the principal household earner} &
\multicolumn{3}{l}{{\scriptsize (ref: Blue collar workers)}}\\
{\scriptsize \quad self-employed persons} & ${\scriptsize \gamma}_{8}$ &
{\scriptsize -0.122} & {\scriptsize 0.771}\\
{\scriptsize \quad white collar workers} & ${\scriptsize \gamma}_{9}$ &
{\scriptsize -3.733} & {\scriptsize
$<$%
0.001}\\
{\scriptsize \quad retirees} & ${\scriptsize \gamma}_{10}$ &
{\scriptsize -5.261} & {\scriptsize
$<$%
0.001}\\
{\scriptsize \quad no activity} & ${\scriptsize \gamma}_{11}$ &
{\scriptsize -1.910} & {\scriptsize 0.004}\\
{\scriptsize Level of Education of the principal household earner } &
\multicolumn{3}{l}{{\scriptsize (ref: BAC and higher degree)}}\\
{\scriptsize \quad student} & ${\scriptsize \gamma}_{12}$ &
{\scriptsize 5.901} & {\scriptsize
$<$%
0.001}\\
{\scriptsize \quad no or weak diploma} & ${\scriptsize \gamma}_{13}$ &
{\scriptsize -1.281} & {\scriptsize
$<$%
0.001}\\\hline
\end{tabular}%
\end{table}%
%

\begin{table}[h] \centering
\caption
{The different steps performed in testing the socioeconomic part of our model. For each step, the null hypothesis tested and the p-value resulting from the appropriate F-test are shown. All tests are performed conditionally on the results of the previous tests (Pval)\label
{Test_socio}}%
\begin{tabular}
[c]{ll}\hline\hline
{\scriptsize Null hypothesis} & {\scriptsize Pval}\\\hline
{\scriptsize H1 : }$\gamma_{8}=0$ & {\scriptsize 0.771}\\
{\scriptsize H2 : }$\gamma_{9}=\gamma_{10}$ & {\scriptsize 0.030}\\
{\scriptsize H3 : }$\gamma_{9}=\gamma_{11}$ & {\scriptsize 0.018}\\
{\scriptsize H4 : }$\gamma_{10}=\gamma_{11}$ & {\scriptsize
$<$%
0.001}\\
{\scriptsize H5 : }$\gamma_{4}=\gamma_{5}=\gamma_{6}=\gamma_{7}$ &
{\scriptsize
$<$%
0.001}\\
{\scriptsize H6 : \ a : }$\gamma_{4}=\gamma_{5}$ & {\scriptsize
$<$%
0.001}\\
{\scriptsize \ \ \ \ \ \ \ \ b : }$\gamma_{4}=\gamma_{6}$ & {\scriptsize
$<$%
0.001}\\
{\scriptsize \ \ \ \ \ \ \ \ c : }$\gamma_{4}=\gamma_{7}$ &
{\scriptsize 0.881}\\
{\scriptsize \ \ \ \ \ \ \ \ d : }$\gamma_{5}=\gamma_{6}$ & {\scriptsize
$<$%
0.001}\\
{\scriptsize \ \ \ \ \ \ \ \ e : }$\gamma_{5}=\gamma_{7}$ & {\scriptsize
$<$%
0.001}\\
{\scriptsize \ \ \ \ \ \ \ \ f : }$\gamma_{6}=\gamma_{7}$ &
{\scriptsize 0.0103}\\
{\scriptsize H7 : }$\gamma_{12}=\gamma_{13}$ & {\scriptsize
$<$%
0.001}\\
{\scriptsize H8 : }$\gamma_{1}=\gamma_{2}=\gamma_{3}$ & {\scriptsize
$<$%
0.001}\\
{\scriptsize H9 : a : }$\gamma_{1}=\gamma_{2}$ & {\scriptsize
$<$%
0.001}\\
{\scriptsize \ \ \ \ \ \ \ \ b : }$\gamma_{1}=\gamma_{3}$ & {\scriptsize
$<$%
0.001}\\
{\scriptsize \ \ \ \ \ \ \ \ c : }$\gamma_{2}=\gamma_{3}$ & {\scriptsize
$<$%
0.001}\\\hline
\end{tabular}%
\end{table}%
%

\begin{table}[h] \centering
\caption
{Restricted maximum likelihood estimates (REML) for all age and socioeconomic variables of the reduced final model with all variance components and their standard errors (s.e)\label
{Esti_final_model}}%
%

\begin{tabular}
[c]{llll}\hline\hline
{\scriptsize Effect} & {\scriptsize Parameter} & {\scriptsize REML} &
{\scriptsize Pval}\\
{\scriptsize Income } & \multicolumn{3}{l}{{\scriptsize (ref: Mean sup)}}\\
{\scriptsize \quad Well to do} & ${\scriptsize \gamma}_{1}$ &
{\scriptsize 6.108} & {\scriptsize
$<$%
0.001}\\
{\scriptsize \quad Mean inf} & ${\scriptsize \gamma}_{2}$ &
{\scriptsize 2.760} & {\scriptsize
$<$%
0.001}\\
{\scriptsize \quad Modest} & ${\scriptsize \gamma}_{3}$ & {\scriptsize -1.915}
& {\scriptsize
$<$%
0.001}\\
{\scriptsize Region of residence} & \multicolumn{3}{l}{{\scriptsize (ref: Non
coastal regions)}}\\
{\scriptsize \quad Paris and North, Brittany, Vendee coast} &
${\scriptsize \gamma}_{4}{\scriptsize =\gamma}_{7}$ & {\scriptsize 0.995} &
{\scriptsize
$<$%
0.001}\\
{\scriptsize \quad South west coast} & ${\scriptsize \gamma}_{5}$ &
{\scriptsize 5.156} & {\scriptsize
$<$%
0.001}\\
{\scriptsize \quad Mediterranean coast } & ${\scriptsize \gamma}_{6}$ &
{\scriptsize 2.250} & {\scriptsize
$<$%
0.001}\\
{\scriptsize Occupation category of the principal household earner} &
\multicolumn{3}{l}{{\scriptsize (ref: Blue collar workers and self employed
persons)}}\\
{\scriptsize \quad white collar workers} & ${\scriptsize \gamma}_{9}$ &
{\scriptsize -3.745} & {\scriptsize
$<$%
0.001}\\
{\scriptsize \quad retirees} & ${\scriptsize \gamma}_{10}$ &
{\scriptsize -5.243} & {\scriptsize
$<$%
0.001}\\
{\scriptsize \quad no activity} & ${\scriptsize \gamma}_{11}$ &
{\scriptsize -1.871} & {\scriptsize 0.005}\\
{\scriptsize Level of education of the principal household earner} &
\multicolumn{3}{l}{{\scriptsize (ref: BAC and higher degree)}}\\
{\scriptsize \quad student} & ${\scriptsize \gamma}_{12}$ &
{\scriptsize 5.879} & {\scriptsize
$<$%
0.001}\\
{\scriptsize \quad no or weak diploma} & ${\scriptsize \gamma}_{13}$ &
{\scriptsize -1.279} & {\scriptsize
$<$%
0.001}\\\hline
&  & {\scriptsize REML} & {\scriptsize s.e}\\
{\scriptsize Variance of the random effect} & ${\scriptsize \sigma}_{u}$ &
{\scriptsize 24.832} & {\scriptsize 6.7316}\\
{\scriptsize Variance-covariance structure } &  &  & \\
{\scriptsize \quad variance} & ${\scriptsize \sigma}^{{\scriptsize 2}}$ &
{\scriptsize 1260705} & {\scriptsize 282309}\\
{\scriptsize \quad correlation} & ${\scriptsize \rho}$ & {\scriptsize -0.22} &
{\scriptsize 0.0434}\\\hline
\end{tabular}%
\end{table}%
%

\end{document}